\documentclass[dvipsnames,format=sigconf]{acmart}
\AtBeginDocument{%
  }

\copyrightyear{2025}
\acmYear{2025}
\setcopyright{cc}
\setcctype{by}
\acmConference[GECCO '25]{Genetic and Evolutionary Computation
Conference}{July 14--18, 2025}{Malaga, Spain}
\acmBooktitle{Genetic and Evolutionary Computation Conference (GECCO '25),
July 14--18, 2025, Malaga, Spain}\acmDOI{10.1145/3712256.3726408}
\acmISBN{979-8-4007-1465-8/2025/07}

\usepackage{cleveref}
\usepackage{nicefrac}

\newcommand{\fid}{\texttt{fid}}
\newcommand{\iid}{\texttt{iid}}

\usepackage{subfig}

\usepackage{algorithm}
\usepackage[noend]{algpseudocode}
\usepackage{algorithmicx}

\algblockdefx{Input}{EndInput}{\textbf{Input: }}{}
\algblockdefx{Output}{EndOutput}{\textbf{Output: }}{}
\algnotext{EndInput}
\algnotext{EndOutput}

\begin{document}

\title[Greedy Restart Schedules]{Greedy Restart Schedules: A Baseline for Dynamic Algorithm Selection on Numerical Black-box Optimization Problems}

\author{Lennart Schäpermeier}
\orcid{0000-0003-3929-7465}
\affiliation{%
  \institution{TU Dresden \& ScaDS.AI Dresden/Leipzig}
  \city{Dresden}
  \country{Germany}
}
\email{lennart.schaepermeier@tu-dresden.de}


\begin{abstract}
In many optimization domains, there are multiple different solvers that contribute to the overall state-of-the-art, each performing better on some, and worse on other types of problem instances.
Meta-algorithmic approaches, such as instance-based algorithm selection, configuration and scheduling, aim to close this gap by extracting the most performance possible from a set of (configurable) optimizers.
In this context, the best performing individual algorithms are often hand-crafted hybrid heuristics which perform many restarts of fast local optimization approaches.
However, data-driven techniques to create optimized restart schedules have not yet been extensively studied.

Here, we present a simple scheduling approach that iteratively selects the algorithm performing best on the distribution of unsolved training problems at time of selection, resulting in a problem-independent solver schedule.
We demonstrate our approach using well-known optimizers from numerical black-box optimization on the BBOB testbed, bridging much of the gap between single and virtual best solver from the original portfolio across various evaluation protocols.
Our greedy restart schedule presents a powerful baseline for more complex dynamic algorithm selection models.
\end{abstract}

\ccsdesc[500]{Mathematics of computing~Continuous optimization}
\ccsdesc[500]{Computing methodologies~Randomized search}
\ccsdesc[500]{Theory of computation~Scheduling algorithms}
\keywords{Heuristic optimization, Continuous optimization, Algorithm scheduling, Dynamic algorithm selection, Benchmarking, Performance analysis}


\maketitle

\section{Introduction}

In heuristic optimization, two different streams of research are constantly at odds with each other: One tries to select, design and configure an algorithm for optimal performance on a specific set of problems, while the other seeks to achieve the best generalist performance for the broadest class of problems.
That both goals cannot be met simultaneously is famously summarized in no-free-lunch theorems which state that there cannot be one universally best optimizer - to improve in one subset of problems, there is a sacrifice of performance in another area \cite{wolpert1997no}.
In the necessary trade-off, a generalist algorithm will always perform somewhat worse than each specialist algorithm in its designated area of expertise.

To close this performance gap, many different methodologies have been developed by different research communities.
\emph{Automated algorithm selection (AAS)}, for example, tries to automate the role of a domain expert with a machine learning model that selects an algorithm for the problem instance at hand \cite{xu2008satzilla,kerschke2019automatedas}.
The AAS model is trained to predict the runtimes of a set of performance-complementary algorithms on the basis of cheap-to-compute numerical instance features, the most promising one being selected to run on the actual optimization problem.
Ideally, this model can efficiently select the best algorithm, only incurring additional costs in the feature computation and a forward pass of the model.
In practice, however, AAS models are costly to train, can often only bridge a small part of the specialist-generalist gap, and suffer from generalization issues when extending beyond their designated training set of instances \cite{kerschke2019automatedbbo,dietrich2024impact,cenikj2024cross}.

An earlier methodological development was the adoption of \emph{principled scheduling and restarting techniques} that incorporate changing population sizes (e.g., IPOP-CMA-ES \cite{auger2005restart} and BIPOP-CMA-ES \cite{hansen2009benchmarking}), or specialized pre-solvers that quickly solve certain problem classes (e.g., HCMA \cite{loshchilov2013bi} and pre-solvers in SATzilla \cite{xu2008satzilla}).
While giving good and robust off-the-shelf performance, these techniques need to be hand-designed for each optimization domain or require deterministic runtimes, and it is unclear how much performance they leave on the table.

More recently, \emph{dynamic algorithm selection (DAS)}, a research direction that tries to combine the more robust scheduling approaches with the more ambitious AAS models, has emerged \cite{vermetten2020towards,schroder2022switching,vermetten2023switch}.
Instead of computing instance features up-front, a DAS model may start with a default algorithm and use the intermediate solution trajectory for feature computation, determining on-the-fly if a switch to another algorithm (and possibly warmstarting the new algorithm) would be beneficial.
DAS was shown to have a promising performance potential \cite{vermetten2020towards}, but its training is even more involved than AAS models, as the solution trajectory of every involved optimizer influences the further optimization procedure differently.

In this paper, we explore a path between a hand-designed restarting strategy and the complications of training a robust DAS model.
In particular, we introduce a data-driven restart scheduling technique: Data from an extensive benchmark run is used to derive a model-free static solver schedule.
To the author's best knowledge, there exists little literature on this scenario, i.e., scheduling heuristic optimizers with varying running times and only probabilistic success.
The method introduced in this report, called greedy restart schedule (GRS), should be considered as a baseline for any algorithm scheduler in this setting.
As a test case, we consider multiple CMA-ES configurations with different population sizes inspired by IPOP-CMA-ES and a selection of other local search approaches, to solve numerical black-box optimization problems.

The paper is structured as follows: We begin by discussing some methodological background from numerical black-box optimization to meta-algorithmic approaches in \Cref{sec:background}.
Then, we dive into a schematic example and the proposed greedy restart schedule in \Cref{sec:schedule}.
Our method is then demonstrated in an experimental study across \Cref{sec:setup,sec:results}.
Finally, we draw conclusions in \Cref{sec:conclusions}.

\section{Background} \label{sec:background}

In this section, we introduce numerical black-box optimization and common solvers in this domain, as well as a short overview of meta-algorithmic approaches, namely feature-based problem characterization, algorithm selection and scheduling, and performance measures employed in the design and evaluation of the meta-algorithmic models.

\subsection{Numerical Black-box Optimization}

In numerical black-box optimization (BBO), the goal is to minimize (w.l.o.g.) $f: \mathcal X \mapsto \mathbb R$, a function that accepts a vector of $n$ numerical parameters $x \in \mathcal X \subseteq \mathbb R^n$ and returns the objective value of the queried point.
No other assumptions about properties of $f$ can be made, such as convexity, unimodality, smoothness, analytical formulation, derivatives, etc.
Numerical BBO problems can be found in application domains such as hyperparameter optimization \cite{loshchilov2016cma,schneider2022hpo} or engineering problems \cite{long2022learning}.

As practical problems often have long evaluation times, unknown properties and proprietary problem definitions, testbeds including test problems with known optima $f^*$ and problem properties are used to benchmark candidate optimizers, comparing their suitability on a broad range of numerical BBO problems.
The most established of these testbeds is the black-box optimization benchmark (BBOB) \cite{hansen2021coco}, which includes 24 functions (\fid s) with countless problem instances (\iid s) per {\fid} to be used for repetitions, as well as an established evaluation protocol.
The BBOB uses a fixed-target evaluation scheme, where the hitting times in terms of function evaluations for different target values $f^* + \tau$ are recorded, and subsequently presented in \emph{runtime profiles}, which are empirical cumulative distribution functions (ECDFs) of solved targets per function evaluation spent.

The most generally well-performing group of optimizers are evolutionary algorithms based on the covariance matrix adaptation evolution strategy (CMA-ES) \cite{hansen2001completely}.
On a surface level, the CMA-ES evolves a sampling distribution, which it updates once per generation based on the success of $\lambda$ candidates sampled from the current distribution.
CMA-ES and variations thereof are predominant optimizers for numerical BBO, owing to its robustness to noise, invariance against rank-preserving transformations, and rotational invariance.
Other search heuristics such as other evolutionary algorithms \cite{tanabe2019reviewing}, local search approaches \cite{powell1964efficient,byrd1995limited}, Bayesian optimizers \cite{hutter2013evaluation}, and hybrids of various differently configured optimizers \cite{loshchilov2013bi,hansen2009benchmarking} can perform favorably on subsets of problems as well.

\subsection{Performance Assessment}

As performance measure for a problem $P$ and algorithm $A$, the expected running time (ERT) \cite{auger2005performance} is often used.
It is computed based on the average running time $t_{A,P}$ and success rate of $p_{A,P} > 0$ per run:
$$
\text{ERT}(A, P) = \frac {t_{A,P}} {p_{A,P}}
$$
Under the assumption that successful and unsuccessful runs have similar runtime distributions, the ERT (as defined above) gives the expected running time until the problem was solved by the given algorithm.
If $p_{A,P} = 0$, the ERT can be defined as infinite, but in meta-algorithmic studies these values are usually imputed with a high penalty value, e.g., a multiple of the maximum runtime as in the penalized averaged runtime (PAR) family of measures \cite{kerschke2019automatedbbo}.

When evaluating the meta-algorithmic approaches (cf. below), performance is usually measured by some aggregate of the ERTs of the individual optimizers.
While a simple average of the ERT may be appealing, it can be biased towards problems that require an inherently large optimization budget.
Addressing this issue, researchers have developed different normalization techniques for the problem-wise ERTs, such as the relative ERT (relERT) \cite{kerschke2019automatedbbo}, which normalizes the ERT by the best known ERT on each problem.
In addition to the relERT, we use the logERT as a performance measure closely related to the area over the curve on the runtime profiles, which has similarly been proposed recently \cite{lopez2024using}.
In the computation of the mean logERT, we compute $\log_{10} \text{ERT}$ before averaging across problems.
Based on the aggregation metric, the single-best solver (SBS) is the algorithm with the best average performance, while the virtual-best solver (VBS), also known as the oracle, comprises the best performance on each problem.

\subsection{Meta-algorithmic Approaches}

In most cases, the same algorithm does not solve all problems the fastest.
To close this gap, meta-algorithmic approaches try to exploit the strengths of multiple algorithms by various means of selecting or combining them.
This includes automated algorithm selection \cite{kerschke2019automatedas,kerschke2019automatedbbo}, which strives to select the best algorithm from a portfolio of solvers, usually based on numerical features describing the optimization landscapes.
In numerical BBO, these are hand-designed exploratory landscape analysis (ELA) features \cite{kerschke2019comprehensive,mersmann2011exploratory}, more recently including end-to-end learned deep ELA features \cite{seiler2025deep}.

Another stream of research dubbed dynamic algorithm selection (DAS) aims to dynamically switch between multiple algorithms during one optimization run \cite{vermetten2020towards,vermetten2023switch,schroder2022switching,guo2024deep}.
Ideally, this can improve performance even more compared to usual algorithm selection as the strengths of multiple algorithms can be combined during the same run, but empirical studies often found it difficult to realize this additional potential \cite{vermetten2023switch,schroder2022switching}.

Finally, there are approaches that don't select one solver but which create solver schedules, e.g., \emph{aspeed} \cite{hoos2015aspeed}.
While timeout-optimal schedules can be derived using techniques such as answer set programming, they require deterministic runtime data (the input is a table of runtimes per instance).
This hinders their application for optimizers that solve instances only with a given probability per run.
Similarly, parallel portfolios of multiple solvers can be utilized, though there are only preliminary studies that don't optimize the underlying portfolio and apply, e.g., a simple round robin approach (where each algorithm is allotted a unit of time in turn) or random restart schemes where a solver from a portfolio is selected at random to run each time \cite{baudivs2014cocopf}.
In practice, due to the dynamic interaction between different optimizers, a parallel portfolio presents more challenges to implement than a restart schedule.

\section{Greedy Restart Schedule} \label{sec:schedule}

We begin with a conceptual example based on an artificial, minimal dataset to illustrate the merits of the scheduling idea.
Afterwards, we present our greedy scheduling heuristic.

\subsection{Schematic Example} \label{sec:motivation}

Let us consider a minimal example where we have two algorithms (A1, A2) that are run on a uniform problem distribution containing problems P1 and P2.
A1 and A2 both require exactly 10 function evaluations each run regardless of problem, but differ in their success rates and subsequently the resulting ERTs: A1 is strong on P1 (20\% success rate) and weak on P2 (5\% success rate), while the success rates for A2 are reversed, cf.~\Cref{tab:minimal_example}. In this example, the VBS has a performance of 50 on each problem, while both A1 and A2 can be considered as the SBS with a mean ERT of 125, respectively.

Let's assume we break ties in favor of A1 and it is run first (due to the symmetry, both choices are equally valid).
A1 is more successful on problem P1, resulting in a \emph{distribution of unsolved problems} after the run that is slightly favoring P2: While 80\% of P1 are left, 95\% of P2 are still unsolved.
As a result, A2 is the more promising algorithm to solve the remaining problem instances, and after a run of A2, we again reach a uniform distribution where now 76\% of both problem classes are unsolved.

\begin{figure*}[t]
    \centering
    \includegraphics[width=\linewidth]{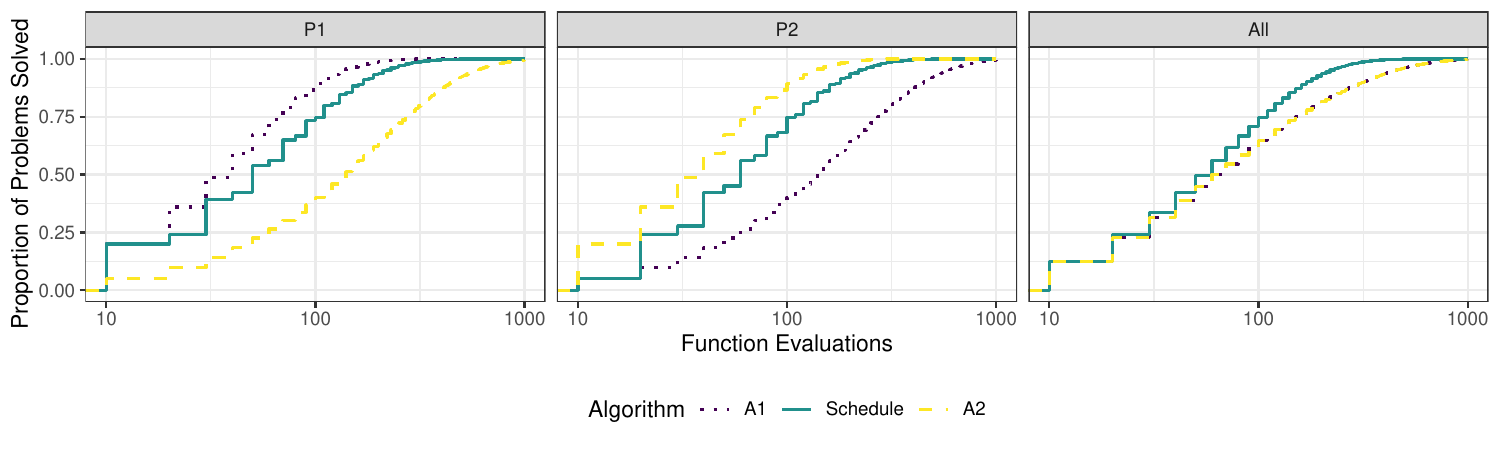}
    \caption{Runtime profiles depicting the proportion of problems solved for A1, A2, and a schedule alternating A1 and A2, for problems P1, P2, and the uniform problem distribution (All). While the simple schedule in this example cannot beat any of its component algorithms on the individual problems, it can close a significant gap to the best respective solver. On the overall problem distribution, it outperforms A1 and A2 (which perform identically).}
    \label{fig:minimal_example}
\end{figure*}

Repeating this schedule ad infinitum, we get the results presented in \Cref{fig:minimal_example}: The schedule is a good second-best on each individual problem, but it outperforms both A1 and A2 considering the whole problem set where P1 and P2 are equally likely to occur.
Further, it achieves an ERT of $75$ on P1 and $81.25$ on P2, resulting in a mean ERT of $78.125$, closing \emph{62.5\% of the SBS-VBS gap}.

This examples goes to show that significant performance gains can be achieved in scheduling optimizers which depend on many, individual restarts to solve optimization problems.
While illustrative, real optimization scenarios contain heterogeneous solvers with widely varying running times and success rates on significantly larger problem sets.
In the following, we introduce a greedy restart scheduling approach that can utilize these larger datasets to create a static solver schedule automatically.

\begin{table}[tb]
    \caption{Expected running times for minimal schematic example from \Cref{sec:motivation}. Best values highlighted in bold.}
    \centering
    \begin{tabular}{c|c|c|c}
         & P1 & P2 & All \\ \hline
        A1 & $\nicefrac{10}{20\%} = \mathbf{50}$ & $\nicefrac{10}{5\%} = 200$ & $125$ \\ \hline
        A2 & $\nicefrac{10}{5\%} = 200$ & $\nicefrac{10}{20\%} = \mathbf{50}$ & $125$ \\ \hline
        Schedule & $75$ & $81.25$ & $\mathbf{78.125}$
    \end{tabular}
    \label{tab:minimal_example}
\end{table}

\subsection{Greedy Scheduler}

Generally, when we evaluate an optimizer using the ERT, we implicitly assume that the underlying algorithm is restarted upon failure to solve the problem to the specified target precision.
Here, we try to utilize this assumption to iteratively select an algorithm that does not solve the initial distribution of problems best, but rather solves the \emph{distribution of unsolved problems at time of selection} best.
The overall algorithm is summarized in \Cref{alg:schedule}.

\begin{algorithm}[!t]
    \caption{Greedy Restart Schedule}\label{alg:schedule}
    \begin{algorithmic}[1]
        
        \Input {Algorithm portfolio $\mathcal A$, problem set $\mathcal P$, running times~$\hat t_{A,P}$ and success rates $\hat p_{A, P} \ \forall A \in \mathcal A, P \in \mathcal P$}
        \EndInput
        
        \Output{Algorithm schedule}
        \EndOutput
        
        \Procedure{GreedyRestartSchedule}{$\mathcal A, \mathcal P, \hat t, \hat p$}
            
            \State $\text{schedule} \gets \text{[]}$
            \State $\pi_P \gets \frac 1 {\|\mathcal P\|} \forall P \in \mathcal P$
            \Comment{(Prior) problem distribution}

            \While{\emph{termination criteria not reached}}
                \State $A^* \gets \arg \max_{A \in \mathcal A} \sum_{P \in \mathcal P} \pi_P \frac {\hat p_{A,P}} {\hat t_{A,P}}$
                \State schedule.append($A^*$)
                \State $\pi_P \gets \pi_P (1 - p_{A^*,P}) \quad \forall P \in \mathcal P$
                \State $\pi_P \gets \frac {\pi_P} {\sum_{P \in \mathcal P} \pi_P} \quad \forall P \in \mathcal P$
            \EndWhile
            
            \State \textbf{return} schedule
    \EndProcedure
    \end{algorithmic}
\end{algorithm}

More precisely, we assume that we have a fixed algorithm portfolio $\mathcal A$, and a (training) set of problem instances $\mathcal P$.
To build the greedy restart schedule, we also require the average runtimes of each algorithm $A \in \mathcal A$ on each problem $P \in \mathcal P$, denoted $t_{A,P}$, as well as the respective success rate per run, denoted $p_{A,P}$.
We approximate those by performing a fixed number of repetitions per algorithm and problem, yielding $\hat t_{A,P}$ as the observed mean time per run and $\hat p_{A,P}$ as the observed proportion of successful runs.
We also initialize a prior distribution across our (training) set of problems.
Here, we assume that all problems from $\mathcal P$ are equally likely, so we use a uniform distribution, i.e., $\pi_P = \frac 1 {\| \mathcal P \|} \forall P \in \mathcal P$.
The prior distribution can be modified to reflect a focus on certain problems, giving some control to model the problem distribution to a user's application at hand.

For the rule by which the next algorithm is selected, there are multiple potential candidates.
Possible options include selecting the algorithm that has the best average ERT (or one of its variants) across problems, though this is likely to favor algorithms that don't have many failures but may require longer running times.
We choose an alternative measure, which is the expected proportion of problems solved per spent evaluation according to the current problem distribution:
$$
A^* = {\arg \max}_{A \in \mathcal A} {\sum_{P \in \mathcal P} \pi_P \frac {\hat p_{A,P}} {\hat t_{A,P}}}
$$
At the beginning, this is more likely to favor specialist algorithms that solve some problem category quickly, and have short runtimes on other problems.
It also does not require an imputation strategy for cases where $\hat p_{A,P} = 0$.
Finally, we update the distribution of unsolved problems using the success rates of the selected algorithm.
By changing the problem distribution, this selection criterion successively focuses on the harder, rarely solved problems and guides the selection of different algorithms over time.

This procedure runs until some termination criterion is met, which can for example be some maximum number of restarts, or when all training problems have been solved with 100\% probability (in which case normalizing $\pi_P$ fails).
The output of our approach is a \emph{static algorithm restart schedule} that does not require domain-specific features or the training of a machine learning model.
To demonstrate the potential of such static schedules, we'll present an experimental study on numerical BBO problems next.



\section{Experimental Setup} \label{sec:setup}

In the following, we introduce problem set, algorithm portfolio and evaluation protocols with which we assess our restart scheduling approach.

\subsection{Problem Set}

As the problem set, we consider test problems from the previously introduced BBOB testbed \cite{hansen2021coco}, namely, all functions $\fid \in \{1,2,\dots,24\}$ with dimensions $d \in \{2,3,5,10\}$.
For training, the $500$ problem instances $\iid \in \{101,102,\dots,600\}$ are chosen, avoiding conflicts with the instances usually chosen when comparing to other algorithms using the regular COCO protocol \cite{hansen2021coco}.

Performance is evaluated based on the fixed-target logERT and relERT measures, where the ERTs of failed runs are imputed $10^7 d$ as a high penalty value corresponding to the cutoff in the runtime profiles.
We consider all default targets of the BBOB problems, i.e., the 51 targets $\tau \in \{10^{2},10^{1.8},\dots,10^{-8}\}$, leading us to create an anytime optimization schedule.
The different problem instances of each function are regarded as equivalent, independent restarts.
Thus, a problem in the context of the scheduler is specified by the tuple $\{\fid,d,\tau\}$, with repetitions given by the 500 instances, resulting in $24 \times 51 = 1\,224$ problems per search space dimension, mirroring the runtime profiles used in COCOs evaluation procedure.

\subsection{Algorithm Portfolio}

The algorithm portfolio to solve the problems is given from two different sources.
On the one hand, we use several configurations of the CMA-ES algorithm.
These are the \texttt{default} configuration of CMA-ES, which uses a population size of $\lambda =  4 + \lfloor 3 \ln d \rfloor$ and variants with doubling population size, up to a 64-fold population size, dubbed \texttt{2L}, \texttt{4L}, \texttt{8L}, \texttt{16L}, \texttt{32L}, and \texttt{64L}, respectively.
These correspond to the configurations which would be successively evaluated by the IPOP-CMA-ES.
Otherwise, we use default parameters and the default termination criteria implemented in the \texttt{modCMA} Python package \cite{denobel2021}.
On the other hand, several local optimizers are chosen from the \texttt{scipy} Python package \cite{virtanen2020scipy}, namely \texttt{SLSQP} \cite{kraft1988software}, \texttt{Powell} \cite{powell1964efficient}, and \texttt{L-BFGS-B} \cite{byrd1995limited}.
These are particularly suited to the unimodal functions with moderate and high conditioning present in the BBOB problem set and converge quickly to a local optimum, presenting a good complement to the CMA-ES variants which are more robust on the (highly) multimodal functions \cite{varelas2019benchmarking}.

\begin{figure*}[t]
    \centering
    \includegraphics[width=0.4\linewidth]{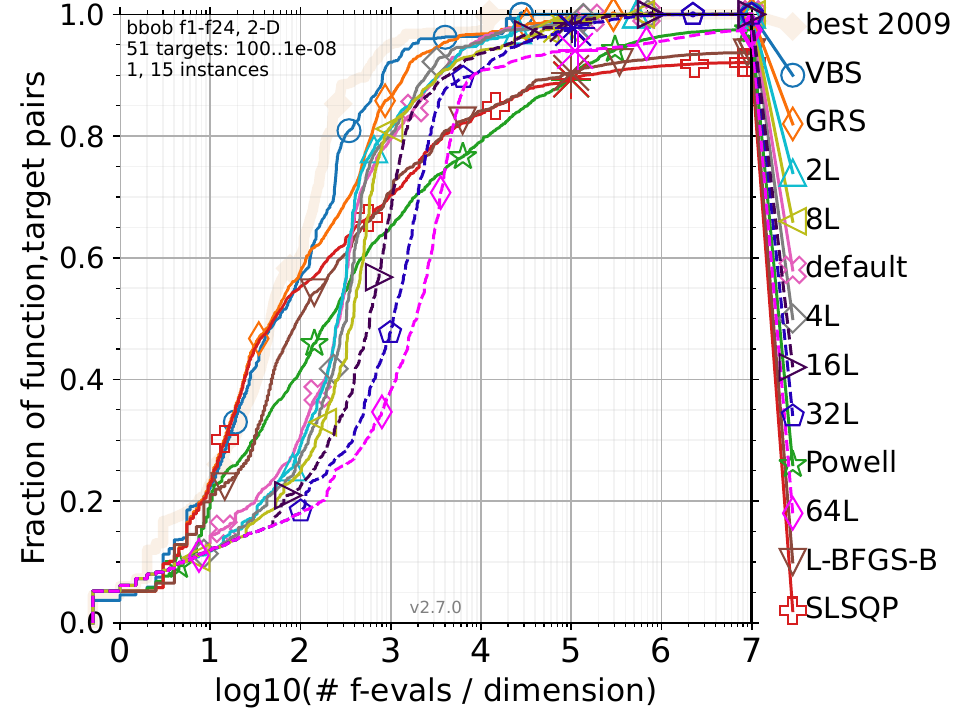}
    \includegraphics[width=0.4\linewidth]{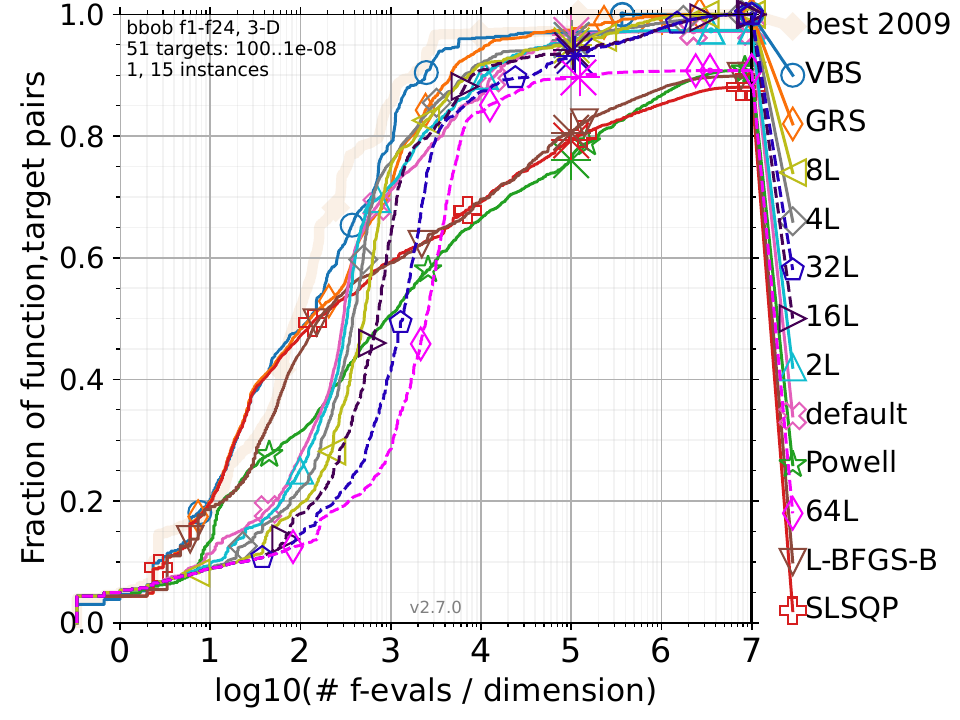}
    \includegraphics[width=0.4\linewidth]{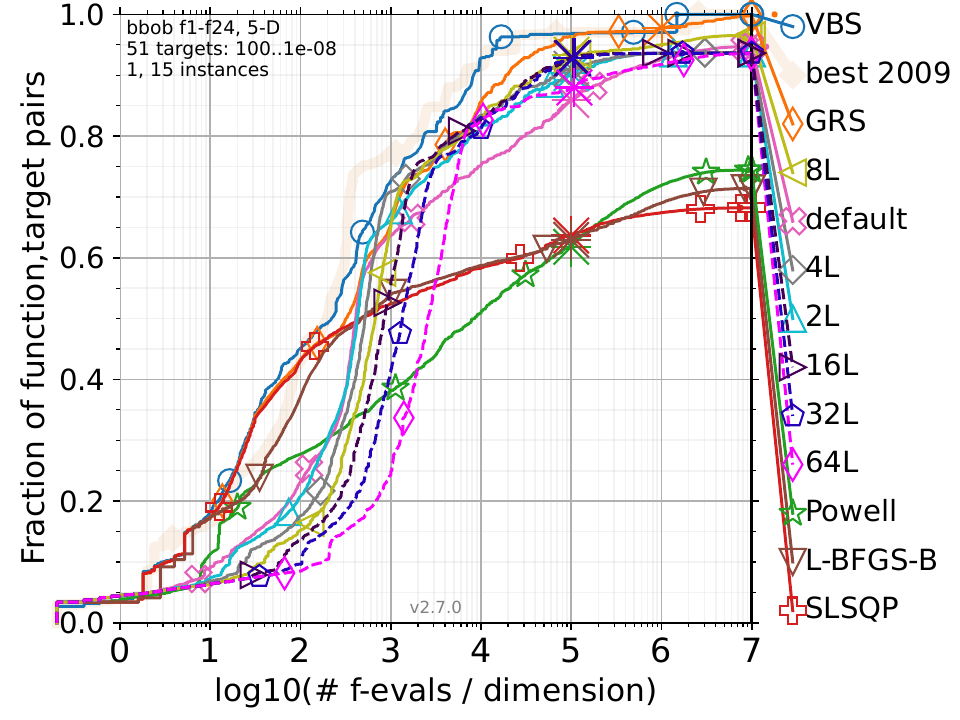}
    \includegraphics[width=0.4\linewidth]{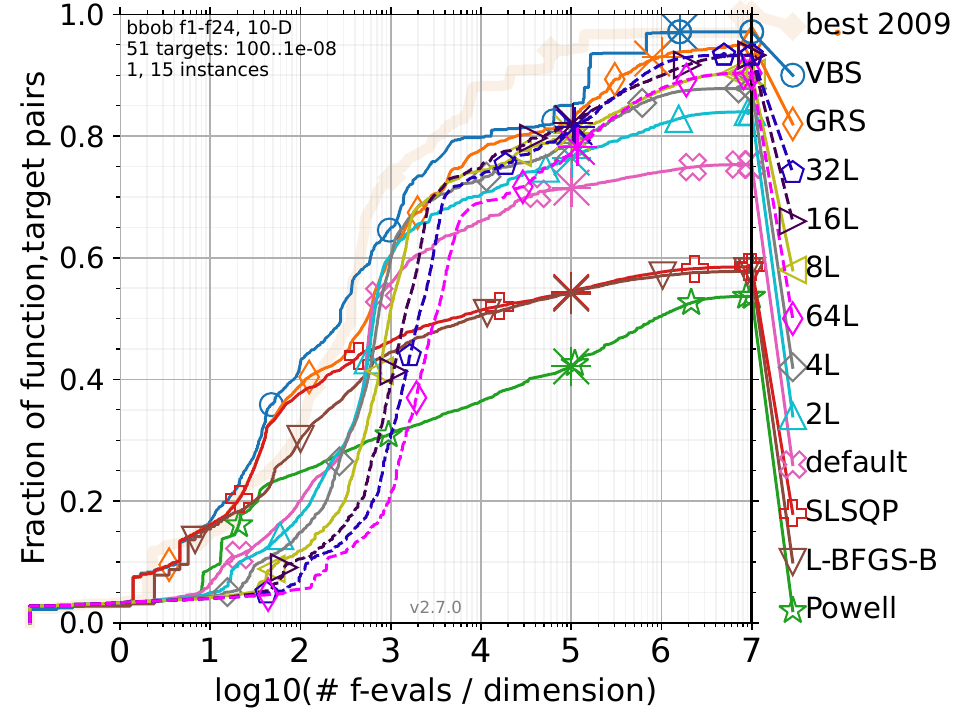}
    \caption{Runtime profiles for $d \in \{2,3,5,10\}$ for all algorithms, the synthetic oracle algorithm (VBS), and the greedy restart schedule (GRS). The schedule beats all individual algorithms overall and approaches the performance of the oracle.}
    \label{fig:grs-performance}
\end{figure*}

\subsection{Evaluation Protocols}

To demonstrate the potential of the scheduling approach, we perform the initial experiments using an instance-based resampling strategy.
This means that the training and testing distribution of problems will be similar (containing all \fid s), though the restriction to a static schedule across all problems provides some protection against overfitting.
To be more precise, after training on \iid s $101-600$, we evaluate on \iid s $1-15$ using the usual COCO protocol, and examine the relERT and logERT on the test set.

We then test our scheduler by resampling using a leave-one-problem-out (LOPO) strategy, where only 23 \fid s are used for training, while evaluating performance on the remaining \fid.
Using LOPO, we can observe the impact of a changing problem distribution on the merits of our approach.
We'll still use the training data from \iid s $101-600$ and test on \iid s $1-15$.

In all our experiments, as the number of variables is trivially known beforehand, a separate schedule is created per search space dimension setting.
Further, all algorithms are run for $100\,000 d$ evaluations during the data collection phase, restarting with different initializations until the budget is exhausted.
The greedy restart schedule (GRS) is then finally run for a maximum number of $1\,000\,000 d$ evaluations.
The source code for all experiments as well as the experimental data used in the scheduling is published at \url{https://github.com/schaepermeier/greedy-restart}.

\section{Experimental Results} \label{sec:results}

Here, we present the experimental results.
We start by discussing the overall performance of the portfolio algorithms, as well as the GRS derived using our method.
Then, we'll analyze properties of the resulting schedule.
Finally, we examine the performance of the GRS in contrast to other state-of-the-art heuristics, before examining the generalization capabilites in the context of the BBOB testbed.

\subsection{Overall Performance}

\begin{figure}
    \centering
    \includegraphics[width=\linewidth]{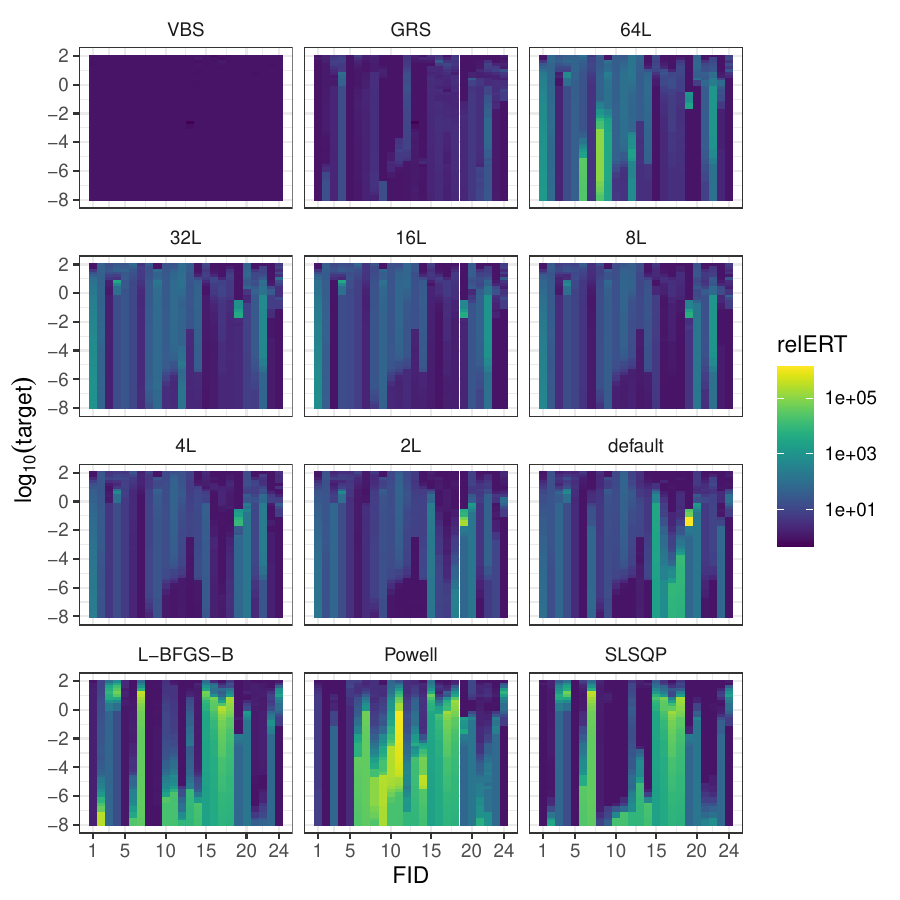}
    \caption{Runtime heatmaps depicting the relERT for all problems with $d=10$. Our GRS provides the best balance of performance across all problems.}
    \label{fig:heatmap-d10}
\end{figure}

The overall runtime profiles for all dimensions are presented in \Cref{fig:grs-performance}.
They include the performance of the individual algorithms, the GRS, as well as a synthetic oracle algorithm (VBS) that consists of the performance data of the algorithm with the best ERT per problem.
We can see that the GRS performs consistently better across the problem distribution than its individual component algorithms, and its performance is always close to the synthetic best algorithm.
\Cref{fig:heatmap-d10} illustrates the performance of the GRS against the other portfolio algorithms for each problem with $d = 10$ using the relERT relative to the portfolio's VBS.
Clearly, the GRS provides the best balance of performance across all individual problems compared to the other algorithms in its portfolio.

In \Cref{tab:mean-erts}, we list the aggregated performances of the GRS for the relERT and logERT metrics, alongside the VBS and SBS from our portfolio.
The SBS is calculated separately for each dimension (as it is a known input) and also for each metric, as the different aggregation techniques have slightly different preferences w.r.t. performance distribution.
In all cases, we can see a reduction in runtime using the GRS when compared to the SBS: Depending on the metric chosen, the overall SBS-VBS gap is closed by more than 95\% (relERT) or 63\% (logERT).

\begin{table}[!t]
    \caption{Mean relERT and logERT performances of GRS and the dimension- and metric-specific SBS, as well as VBS. The GRS closes the SBS-VBS gap substantially.}
    \centering

    \subfloat[relERT]{
        \begin{tabular}{crrrr}
            $d$ & GRS & SBS & VBS & Gap closed \\ \hline
            2 & 1.90 & 8.45 & 1.00 & 87.92\% \\
            3 & 2.11 & 12.25 & 1.00 & 90.13\% \\
            5 & 2.56 & 93.81 & 1.00 & 98.32\% \\
            10 & 4.48 & 50.50 & 1.00 & 92.97\% \\ \hline
            all & 2.76 & 41.25 & 1.00 & 95.63\% \\
        \end{tabular}
    }
    \quad
    \subfloat[logERT]{
        \begin{tabular}{crrrr}
            $d$ & GRS & SBS & VBS & Gap closed \\ \hline
            2 & 2.30 & 2.71 & 2.09 & 66.13\% \\
            3 & 2.75 & 3.16 & 2.52 & 64.06\% \\
            5 & 3.26 & 3.73 & 2.96 & 61.04\% \\
            10 & 3.98 & 4.49 & 3.61 & 57.95\% \\ \hline
            all & 3.07 & 3.53 & 2.80 & 63.01\% \\
        \end{tabular}
    }
    \label{tab:mean-erts}
\end{table}

\subsection{Schedule Analysis}

\begin{figure}[t]
    \centering
    \includegraphics[width=\linewidth]{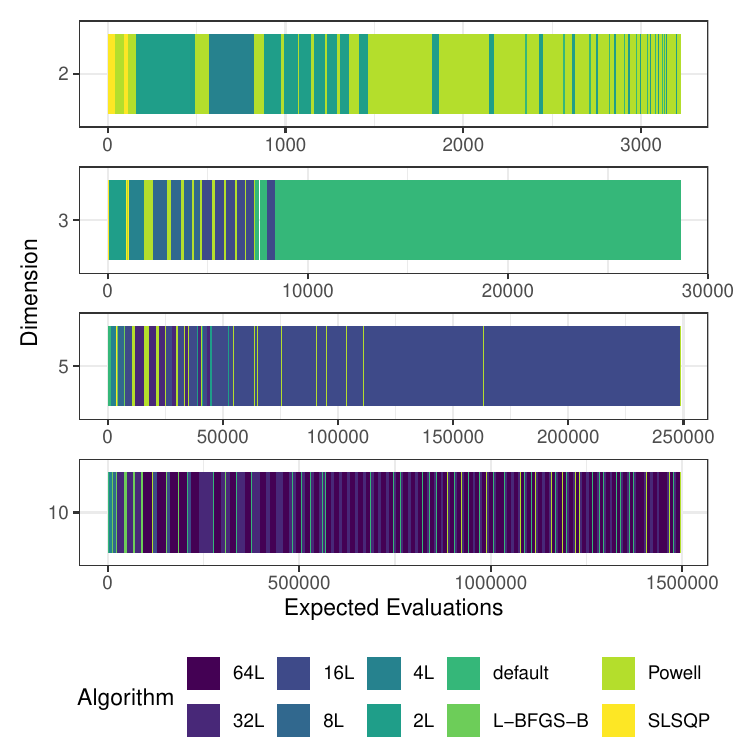}
    \caption{GRS composition for $d \in \{2,3,5,10\}$ showing the first $1\,000$ restarts in each scenario. The algorithms are stacked in order of execution, showing the expected evaluations w.r.t. the BBOB problem set.}
    \label{fig:1k-illustration}
\end{figure}

The composition of the schedules (up to $1\,000$ restarts) is presented in \Cref{fig:1k-illustration}.
Here, we can see that the schedule composition varies quite a bit between the different search space dimensions.
While there is a more diverse mix of solvers run at the beginning, there seems to be some convergence to a particular algorithm or set of algorithms that runs best on the most difficult to solve problems, e.g., Powell for $d=2$ and the default CMA-ES for $d=3$.
It should be noted that only for $d=2$ approximately $100\%$ of all problems are solved, while for $d\in\{5,10\}$ there are even problems that aren't solved by any algorithm in the portfolio (e.g., the more precise targets of \fid~24), which influences this illustration of the results.

\subsection{Comparison with Other Hybrid Heuristics}

\begin{figure*}
    \centering
    \subfloat[Runtime profiles for all \fid s for $d\in \{2,3,5\}$.]{
        \begin{minipage}{\linewidth}
            \centering
            \includegraphics[width=0.32\linewidth]{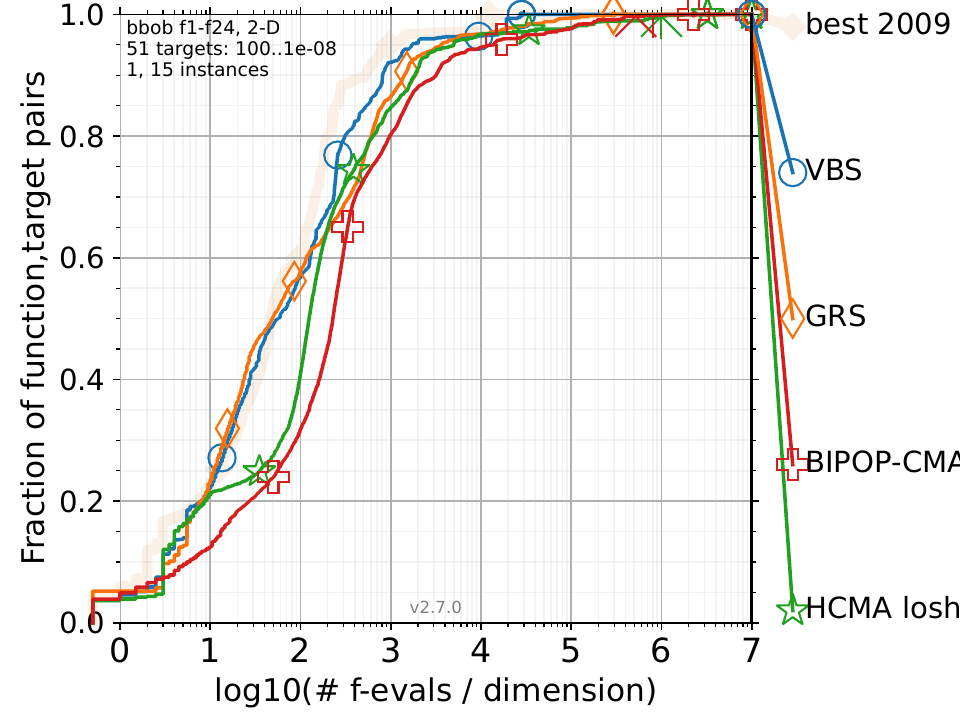}
            \includegraphics[width=0.32\linewidth]{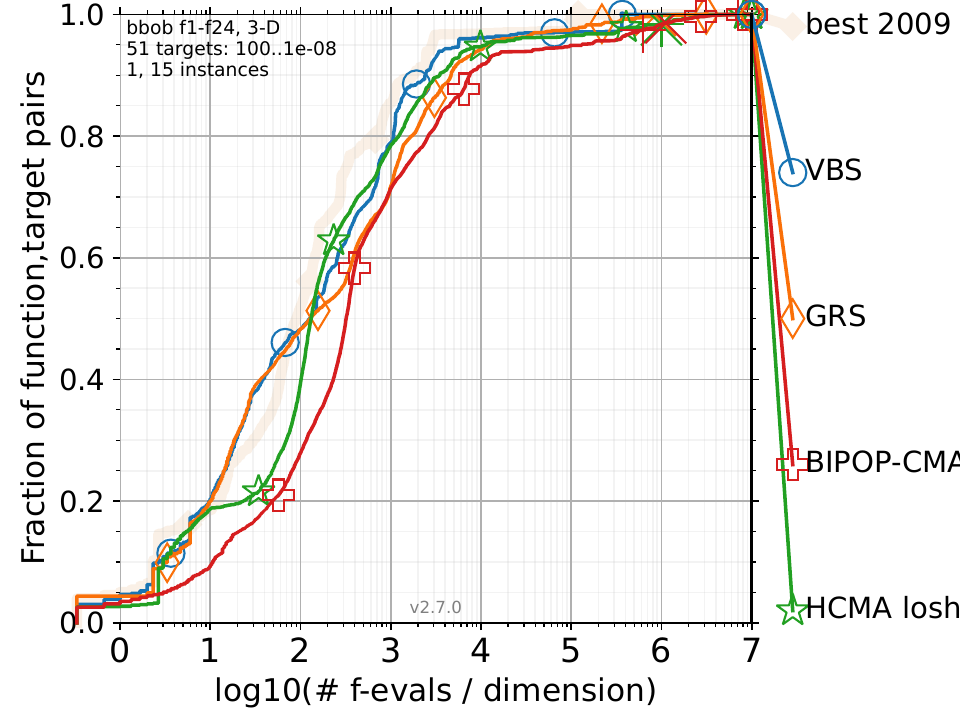}
            \includegraphics[width=0.32\linewidth]{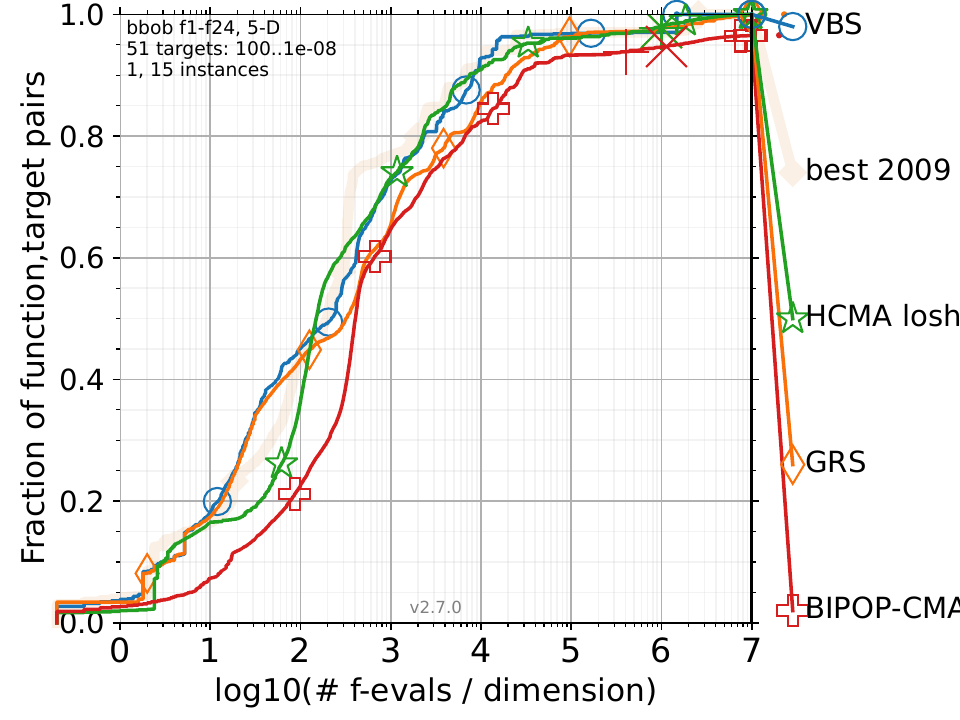}
        \end{minipage}
    }
    
    \subfloat[Runtime profiles per function group and for all \fid s for $d = 10$.]{
        \begin{minipage}{\linewidth}
            \centering
            \includegraphics[width=0.32\linewidth]{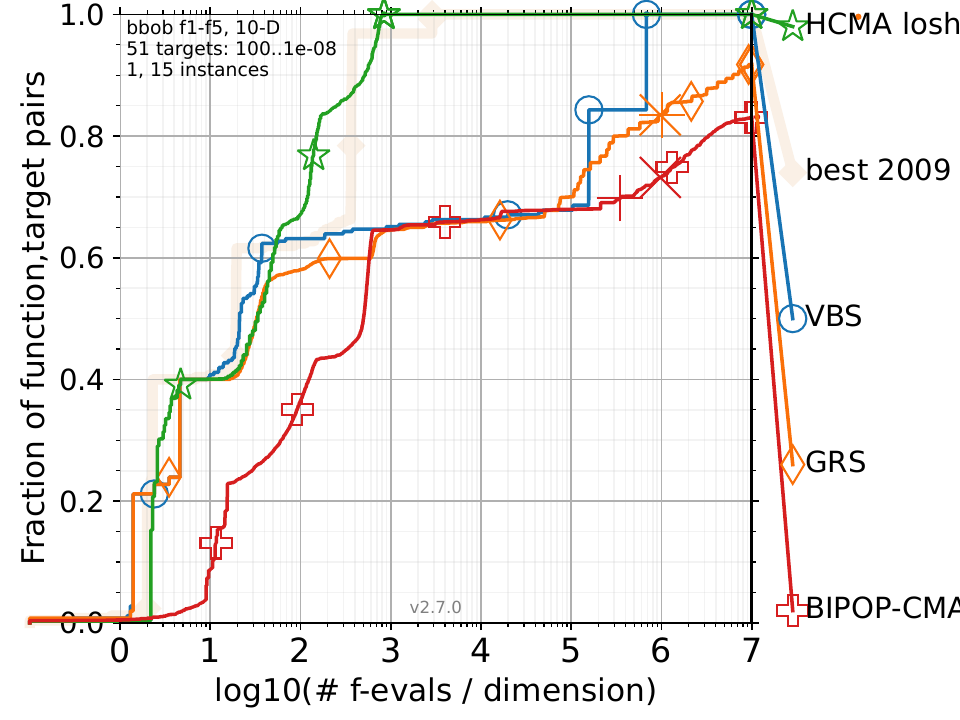}
            \includegraphics[width=0.32\linewidth]{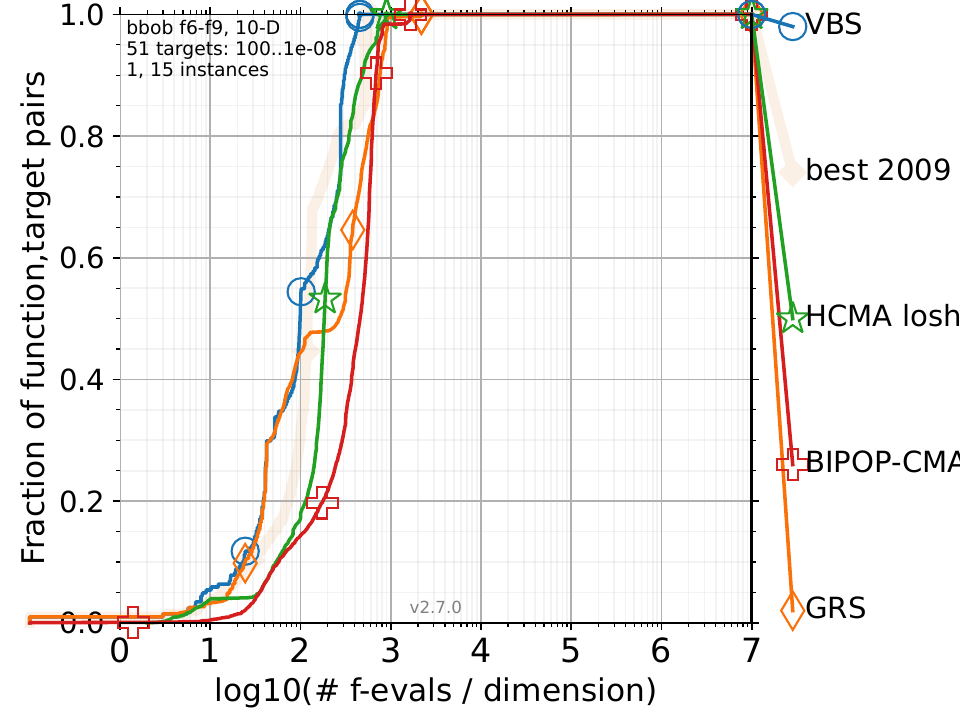}
            \includegraphics[width=0.32\linewidth]{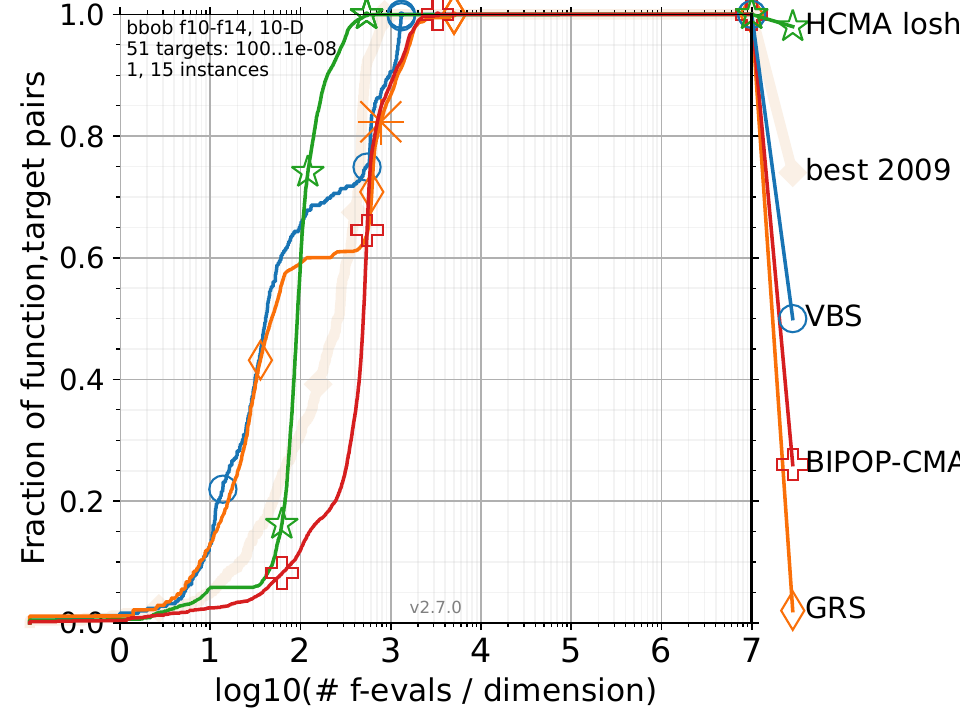}
            \includegraphics[width=0.32\linewidth]{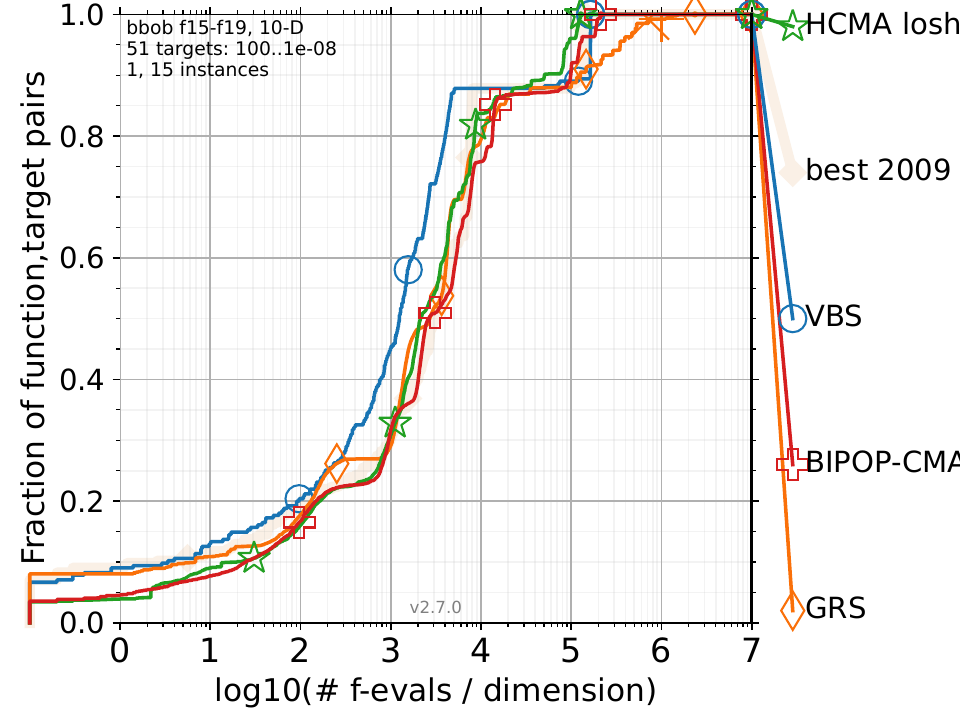}
            \includegraphics[width=0.32\linewidth]{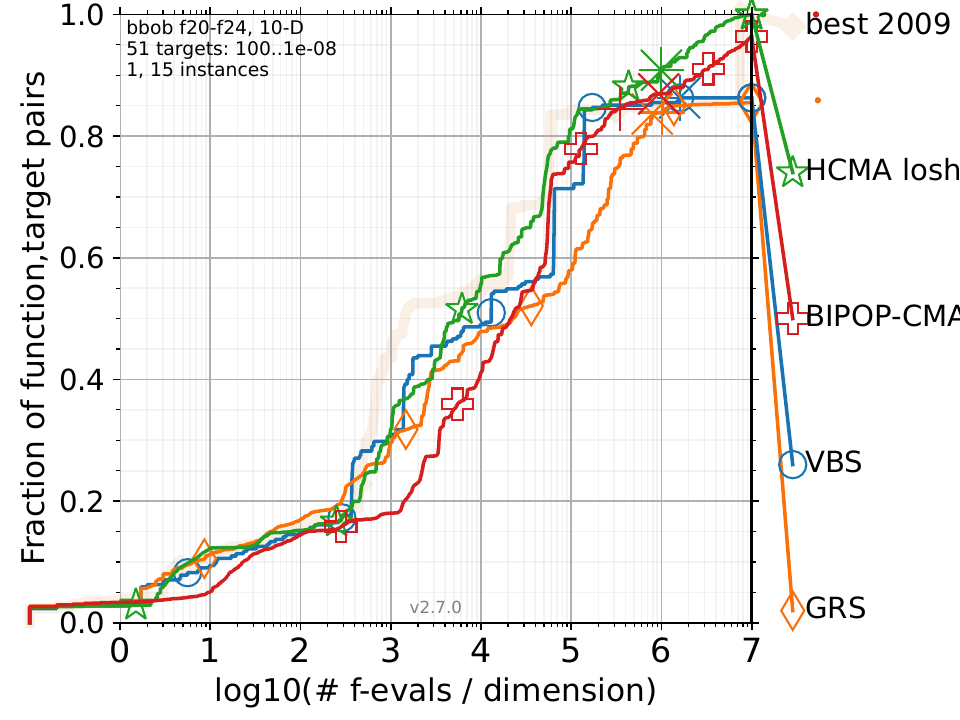}
            \includegraphics[width=0.32\linewidth]{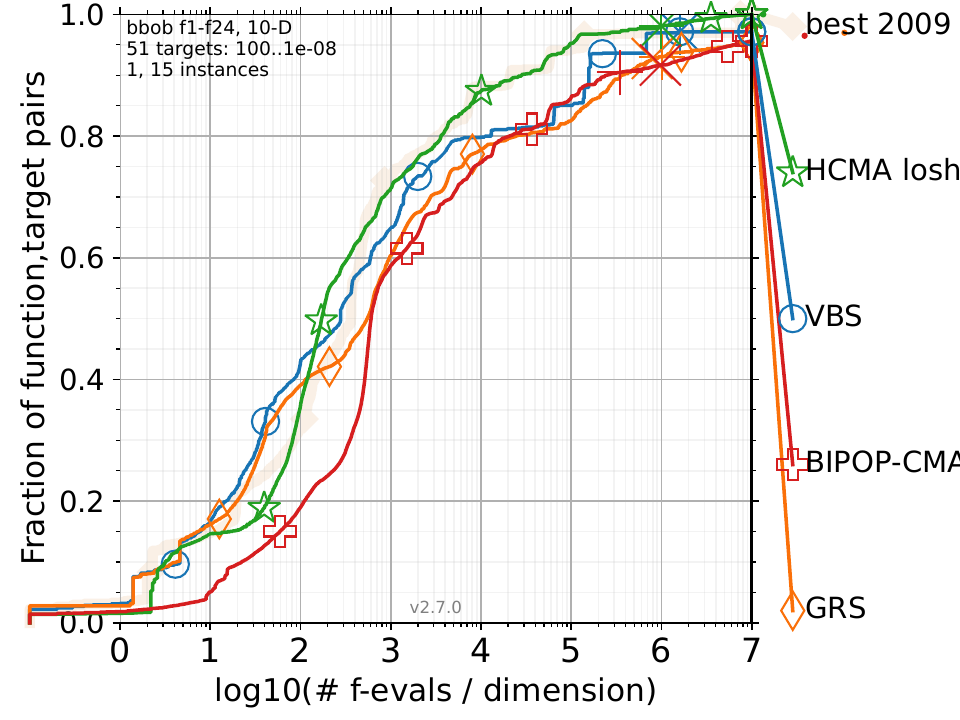}
        \end{minipage}
    }
    \caption{Runtime profiles for GRS, the synthetic oracle algorithm (VBS), and the hybrid optimizers HCMA and BIPOP-CMA-ES (a) across all functions for $d \in \{2,3,5\}$ and (b) per function group and for all functions in $d=10$. GRS beats the other approaches for lower dimensions, and often has a performance advantage on low budgets compared to the other off-the-shelf solvers.}
    \label{fig:hcma-performance}
\end{figure*}

\Cref{fig:hcma-performance} depicts runtime profiles for GRS compared to the hybrid restart heuristics HCMA \cite{loshchilov2013bi} and BIPOP-CMA-ES \cite{hansen2009benchmarking}, which are state-of-the-art solvers for off-the-shelf numerical BBO.
On $d\in\{2,3\}$, GRS presents the best overall results, while still performing well for $d \in \{5,10\}$.
In all scenarios, it has a performance advantage on lower budgets.
The biggest gap is created by HCMA on separable problems, as it includes a component algorithm (multivariate STEP \cite{swarzberg1994step}) suited for multimodal separable problems (\fid s 3 and 4), where even the oracle on our portfolio lags behind in performance.

\subsection{Leave-one-problem-out Generalization}

\begin{figure*}[t]
    \centering
    \includegraphics[width=0.32\linewidth]{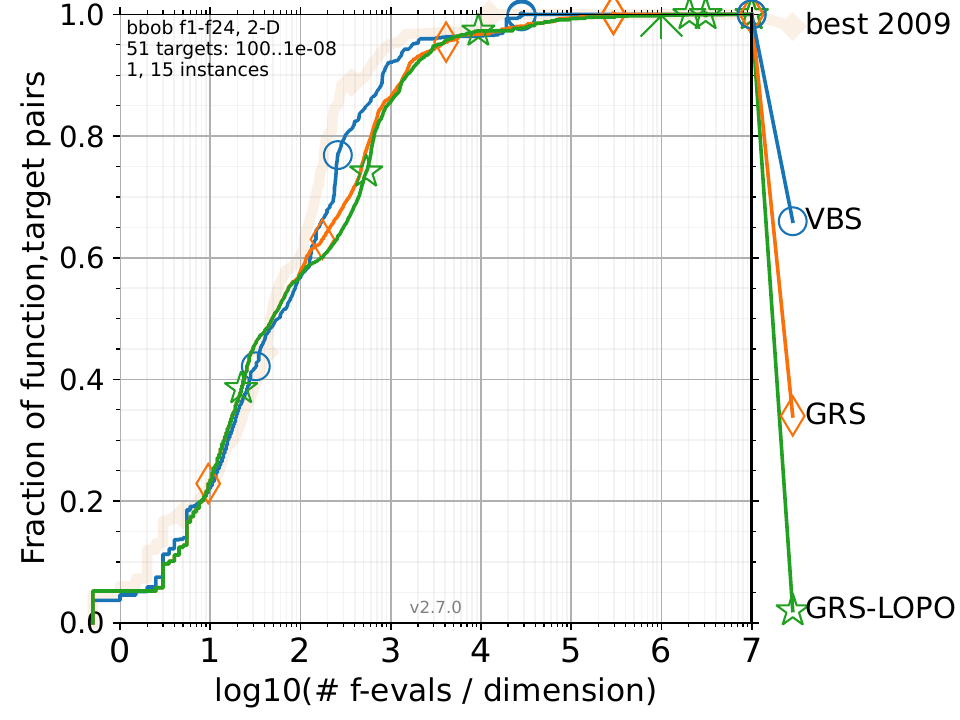}
    \includegraphics[width=0.32\linewidth]{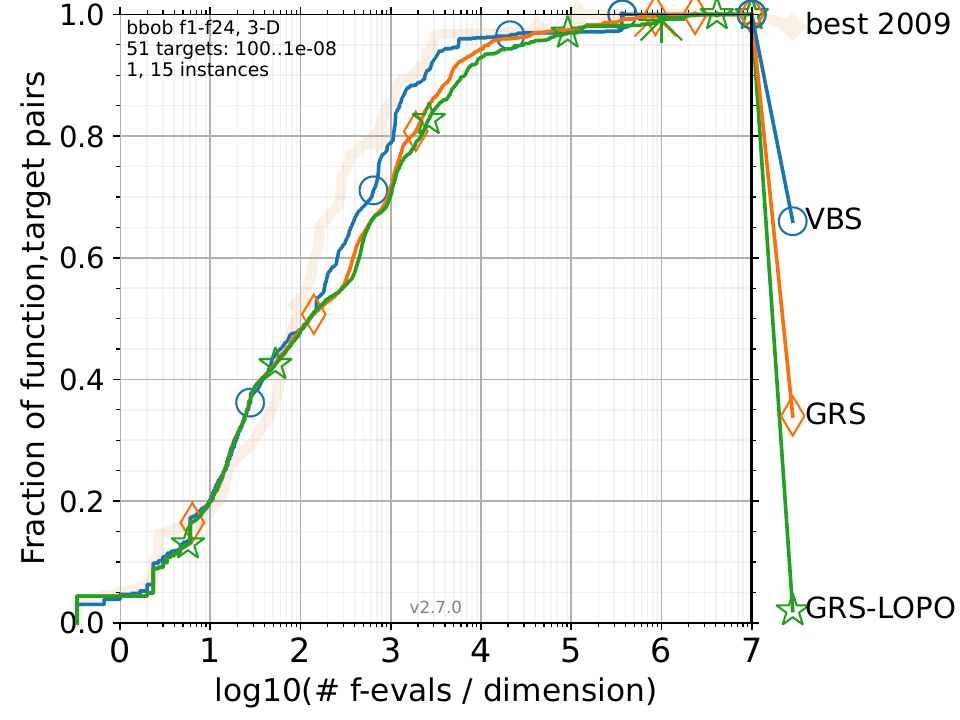} \\
    \includegraphics[width=0.32\linewidth]{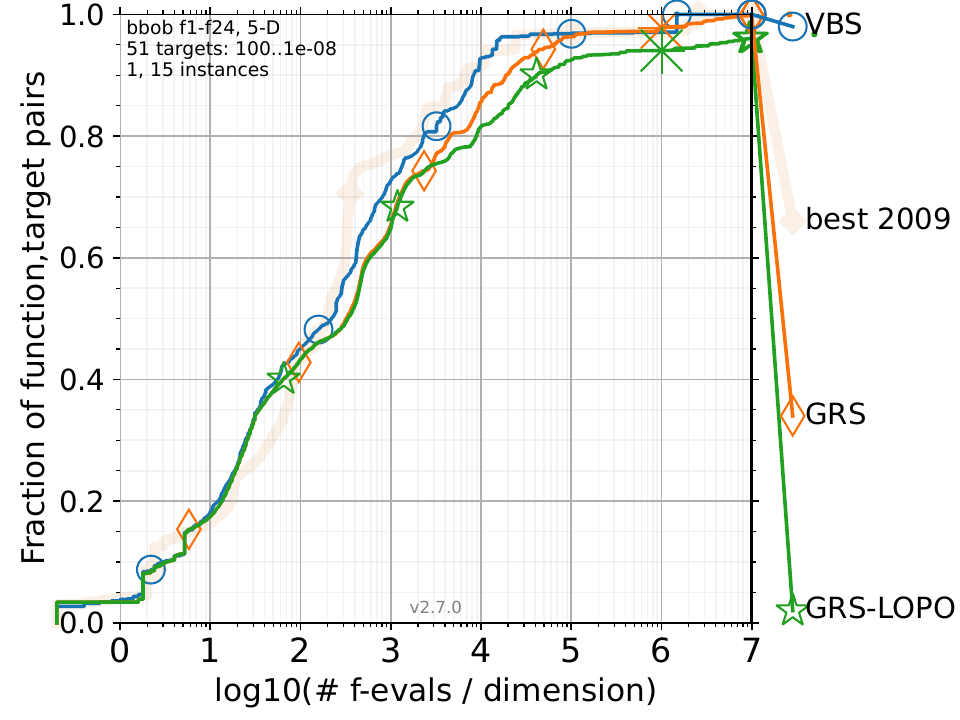}
    \includegraphics[width=0.32\linewidth]{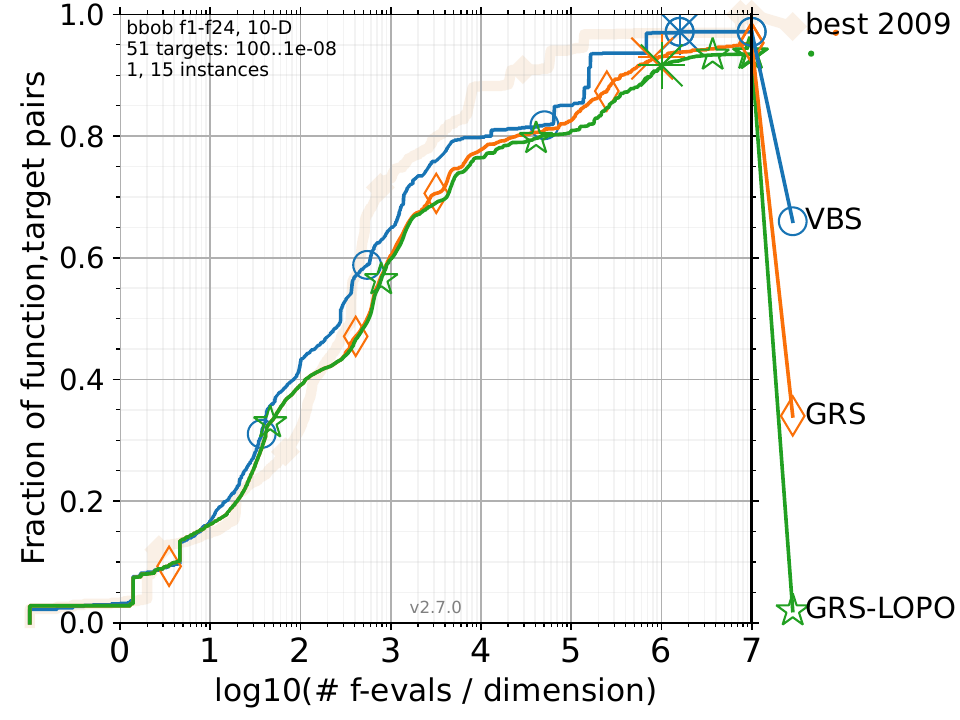}
    \caption{Impact of LOPO resampling on the runtime profiles. There is a noticeable performance gap to the usual GRS starting with $d=5$, but the performance is still close to the oracle.}
    \label{fig:grs-lopo}
\end{figure*}

Here, we present the results of leave-one-problem-out (LOPO) resampling: On the evaluation of each function, the applied schedule is learned from the other 23 \fid s.
The runtime profiles for the BBOB testbed are presented in \Cref{fig:grs-lopo}.
There, we can observe a small performance gap to the unrestricted GRS in higher dimensions, though it is rather small and the overall performance is still close to the oracle.

\begin{figure}
    \centering
    \includegraphics[width=\linewidth]{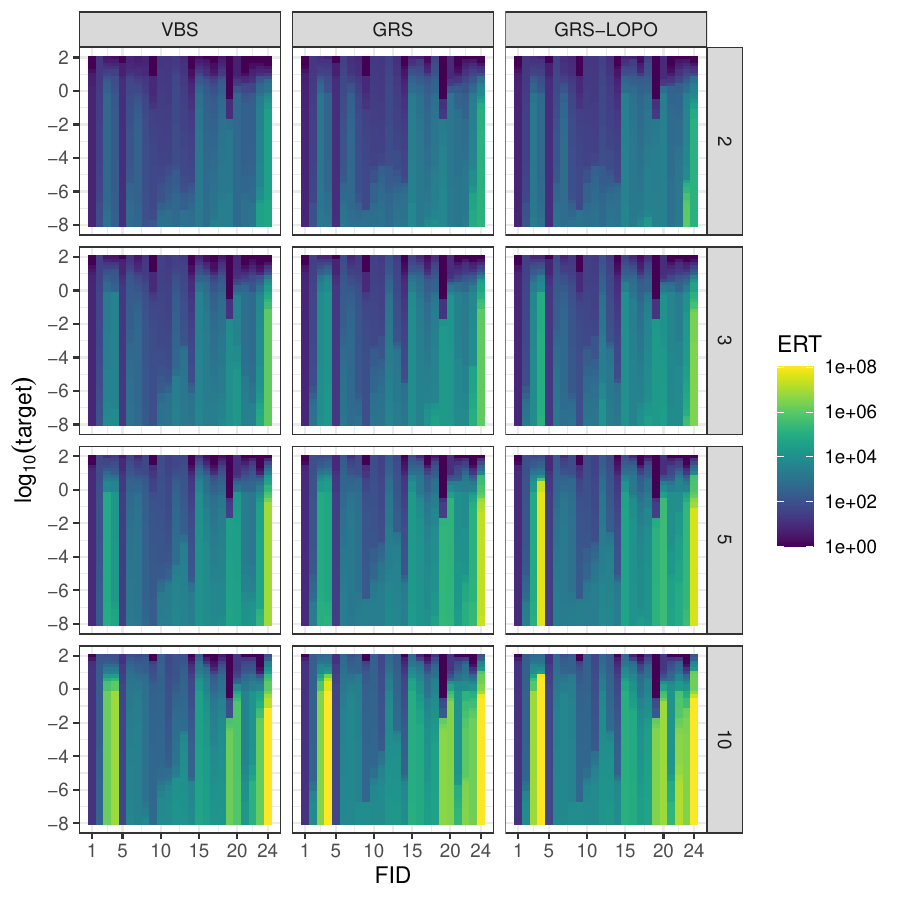}
    \caption{Runtime heatmaps depicting the ERTs of the VBS, GRS with and without LOPO resampling for $d\in\{2,3,5,10\}$. Overall, the performance of GRS and GRS-LOPO is quite similar, with only few noticeable changes, e.g., \fid~4 for $d=5$, where GRS-LOPO fails to solve some targets.}
    \label{fig:heatmap-lopo}
\end{figure}
The performance of the GRS with and without LOPO for all problems is also presented in \Cref{fig:heatmap-lopo}.
The only noticeable differences are a worse performance on \texttt{fid}~4 for $d=5$, which GRS-LOPO did not solve completely, and some higher values in the final function group, e.g., \texttt{fid}~23 for $d=2$.
Due to the LOPO, and thus the different training and testing distributions, we can observe a reduced mean relERT and logERT performance: Overall mean logERT increases to 3.15 (up from 3.07), closing 52\% gap (down from 63\%), while mean relERT increases to 30.28 (up from 2.76), closing 27\% gap (down from 96\%).
While in both cases the schedule still comfortably outperforms the SBS, the large drop in relERT performance is surprising.
This sharp jump can mostly attributed to the worse performance on the already mentioned five-dimensional \fid~4, which has an outsized influence on the relERT, but serves to show that there is still some dependence on the distribution of training problems.

\section{Conclusions} \label{sec:conclusions}

In many optimization domains, deciding which algorithm to choose for a given, but previously unknown, problem is a challenging endeavor.
While many optimizers work best in a restart-oriented fashion, and performance indicators like the expected running time based on repeated restarting are widespread, only little research has so far been dedicated in the literature to extract the most performance using static restart schedules that alternate between solvers.

In this report, we have introduced a greedy restart scheduling (GRS) approach that creates a static restart schedule from a database of previously evaluated optimizers.
The greedy restart schedule has been evaluated on black-box optimization as a test case, using a portfolio of CMA-ES variants derived from the popular IPOP-CMA-ES restart strategy \cite{auger2005restart}, augmented by further local searchers from the \texttt{scipy} Python library \cite{virtanen2020scipy}.
Trained on the BBOB testbed \cite{hansen2021coco}, the static solver schedule is optimized for anytime performance and consistently outperforms the individual algorithms from the portfolio, substantially closing the gap from the single best solver to the virtual best solver by more than 95\% regarding the relative ERT and more than 63\% regarding the logERT.
The performance is robustly improved even under more challenging evaluation protocols such as leave-one-problem-out resampling, though some sensitivity to the underlying problem distribution can be observed.
In a comparison to the state-of-the-art hybrid heuristics HCMA \cite{loshchilov2013bi} and BIPOP-CMA-ES \cite{hansen2009benchmarking}, our GRS performs favorably for lower dimensions and has a consistent performance advantage on lower budgets.
As the algorithm portfolio employed in this study is rather basic, there is still likely performance to be gained by increasing the diversity of optimizers, especially including better solvers for highly multimodal problems.

The approach presented here can be considered as a baseline for any kind of algorithm schedule and dynamic algorithm selection.
Data-driven static schedules also present a basis upon which to build more capable approaches: For example, using the already generated samples to compute landscape features to learn runtime prediction models, giving a potentially even better indication which algorithm to choose for the next restart.

With all these achievements, our scheduling approach presented here is still quite simple.
For example, it is operating in a greedy fashion, so there is likely some potential left in optimizing the static schedules by applying a more sophisticated scheduler.
Our approach is also dependent on reasonable termination criteria of the individual algorithms.
This requirement can be relaxed, if we consider the possibility of selecting not just an algorithm, but also its specific budget for the next run, which presents another avenue for future research.

Finally, we plan to expand the application of data-driven restart schedules to other heuristic optimization domains, such as the inexact solving of traveling salesperson problems where restart heuristics also dominate the state-of-the-art \cite{kerschke2018leveraging,heins2024dancing}.


\bibliographystyle{ACM-Reference-Format}
\bibliography{references}

\end{document}